\documentclass[10pt]{article}

\usepackage{latexsym}

\author{Andrea Surroca}

\title{\Large Sur le nombre de points alg\'ebriques o\`u une fonction analytique transcendante prend des valeurs alg\'ebriques}

\date{\today}

\usepackage[english,francais]{babel}

\newtheorem{theorem}{Theorem}[section]
\newtheorem{lemma}[theorem]{Lemma}
\newtheorem{e-proposition}[theorem]{Proposition}

\newtheorem{e-definition}[theorem]{Definition\rm}

\newtheorem{theoreme}{Th\'eor\`eme}[section]
\newtheorem{lemme}[theoreme]{Lemme}

\newcommand\enteros{\mathbf{Z}}

\newcommand\rationnels{\mathbf{Q}}

\newcommand\reels{\mathbf{R}}

\newcommand\complexes{{\mathbf{C}}}

\newcommand\Qbarre{\overline{\mathbf{Q}}}

\newcommand\Drferme{\overline{D(0,r)}}

\newcommand\Duferme{\overline{D(0,1)}}

\newcommand\DRouvert{D(0,R)}

\newcommand\card{\mathrm{card}}

\newcommand\esp{\hspace{0,2cm}}

\begin{document}

\selectlanguage{francais}

\maketitle

\thispagestyle{empty}


\begin{quote}
\textbf{R\'esum\'e.}
{\small On \'etudie l'ensemble des nombres alg\'ebriques de hauteur
  et de degr\'e born\'es o\`u une fonction analytique transcendante prend des
  valeurs alg\'ebriques.}
\end{quote}


\selectlanguage{english}


\begin{quote}
\textbf{Abstract.}
{\small We study the set of algebraic numbers of bounded height and bounded degree where an analytic transcendental function takes algebraic values.}
\end{quote}


\noindent \textbf{Abridged English version}

\medskip

Given a function, $f,$ analytic and transcendental over $\mathbf{C}[z],$ we are interested on the set, $S_f,$ of algebraic points on which $f$ takes algebraic values. For exemple, for the exponential function, $S_f=\{0\},$ and for $f(z)=2^z, \esp S_f=\rationnels.$ We know, from St\"ackel (\cite{stackel}, \cite{gramain} and \cite{mahler} chap.3), that there exist entire transcendental functions, which take algebraic value at every algebraic point, i.e. for which $S_f = \Qbarre.$ In theorem  \ref{contrex}, we construct an entire  transcendental function $f$ such  that, for all $\alpha \in \Qbarre,$ $f(\alpha) \in \rationnels[\alpha].$ 
   
We classify the countable set of all algebraic numbers by the degree and the height. If $\alpha $ is an algebraic number, we denote its degree by $d(\alpha),$   its Mahler mesure by $M(\alpha)$ and its absolute logarithmic height, $\frac{1}{d(\alpha)}\log M(\alpha),$ by $h(\alpha).$ 

For $R \in \reels \cup \{\infty\}, \esp r \in \reels$ such that  $R>r>0,$ and, $f$ a function analytic over $\DRouvert,$ for every integer $D \geq 1$ and every real $N \geq 0, \esp \Sigma_{D,N}=\Sigma _{D,N}(f,r)$ will denote the set of numbers $\alpha \in \Qbarre \cap  \Drferme$ such that 
$$f(\alpha) \in
\Qbarre ,\esp [\rationnels(\alpha , f(\alpha )):\rationnels] \leq D ,\hspace{0,2cm} h(\alpha) \leq N \esp \mathrm{and} \esp h(f(\alpha)) \leq N,$$
and $\sigma_{D,N}$ will denote its cardinal. So, $S_f$ is the union of all $\Sigma_{D,N}$ for$D\geq 1$ and $N\geq 0.$

\begin{theorem}\label{contrex-e}
Let $\phi$ be a positive function sucht that $\phi(x)/x$ tends towards $0$ when $x \to \infty.$ There  exist an entire  function $f$ transcendental over $\complexes [z ],$ such that
$$f(\alpha) \in \rationnels[\alpha], \esp \forall \alpha \in \Qbarre, $$
and for every integer $D \geq 1,$ there is an infinite number of reals $N \geq 0,$ verifying
$$\card(\Sigma _{D,N}(f,1)) \geq \frac{1}{2}e^{D(D+1)\phi(N)}.$$

\end{theorem}

For fixed $D$ and $N,$  this set is finite. In fact, it is contained in the set
$$E_{D,N}= \{ \alpha \in \Qbarre / \hspace{0,2cm}[\rationnels (\alpha
):\rationnels ] \leq D , \esp h(\alpha )\leq N \}$$
which is finite. We denote by $\epsilon_{D,N}$ its cardinal. 

\medskip 

\begin{lemma}\label{e_{D,N}-e}
For every integer $D \geq 1$ and every real $N \geq 0,$ the cardinal number $\epsilon _{D,N}$ of $E_{D,N}$ is such that
$$e^{D(D+1)(N-1)} < \epsilon _{D,N} \leq e^{D(D+1)(N+1)}.$$

\end{lemma} 

The proof is elementary. The upper bound is easy; a more accurate estimate as $N$ tends towards infinity is known \cite{chern-vaaler}, but is not explicit. To obtain the lower bound, we count the polynomials in  $\enteros[X]$ of bounded height which are 2-Eisenstein over the ring $\enteros_2$ of 2-adic integers; a more accurate upper bound is given in \cite{loher}. It seems that no asymptotic estimate is known (\cite{schmidt} p. 27).

The following result is the main theorem of this paper.

\begin{theorem}{Upper bound for the cardinal number of $\Sigma_{D,N}(f,r).$}\label{propprincipale-e}

Let $R \in \reels \cup \{\infty\}, \esp r \in \reels$ such that $R>r>0,$ and $f$ a function which is analytic over $\DRouvert.$ There is a constant $\gamma >0,$ which depends only on $R, r$ and $f,$ such that, for every integer $D \geq1,$ there is an infinite number of reals $N \geq 0$ for which
$$\card(\Sigma _{D,N}(f,r)) < \gamma D^3N^2.$$

\end{theorem}

Even if we replace the upper bound in Theorem \ref{propprincipale-e}, by any function $\leq \frac{1}{2}e^{D(D+1)\phi(N)},$ Theorem \ref{contrex-e} (resp. \ref{propprincipale-e})  shows that we cannot replace, in Theorem \ref{propprincipale-e} (resp. \ref{contrex-e}), ``there is an infinity of reals $N$'' by ``for all $N$ large enough''. For the particular case $D=1,$ Elkies (\cite{elkies} th.4.) has shown that, for every $\epsilon > 0,$ there exixts a positive constant $A_{\epsilon}$ such that for every real $N\geq 0, \esp \card(\Sigma_{1,N}) \leq A_{\epsilon} e^{\epsilon N}.$

\par\medskip\centerline{\rule{2cm}{0.2mm}}\medskip

\setcounter{section}{0}

\selectlanguage{francais}


\section{Introduction et r\'esultats}

\'Etant donn\'ee une fonction $f$ analytique 
et transcendante sur $\mathbf{C}[z],$ on s'int\'eresse \`a l'ensemble $S_f$ des points alg\'ebriques en lesquels la fonction $f$ prend des valeurs alg\'ebriques. Par exemple, pour la fonction exponentielle, $S_f=\{0\},$ et si on pose $f(z)=2^z,$ alors $S_f=\rationnels.$ On sait, d'apr\`es St\"ackel (cf. \cite{stackel} et aussi \cite{gramain} et \cite{mahler} chap. 3), qu'il existe des fonctions enti\`eres et  transcendantes, prenant des valeurs alg\'ebriques en tous les points alg\'ebriques, c'est-\`a-dire, telles que $S_f = \Qbarre.$ Dans le th\'eor\`eme \ref{contrex}, on construit une fonction $f$ enti\`ere et transcendante v\'erifiant, pour tout $\alpha \in \Qbarre,$ $f(\alpha) \in \rationnels[\alpha].$ 
   
Nous filtrons l'ensemble d\'enombrable des nombres alg\'ebriques par le degr\'e et la hauteur. Si $\alpha $ est un nombre alg\'ebrique, on note $d(\alpha)$ son degr\'e, $M(\alpha)$ sa mesure de Mahler et $h(\alpha)= \frac{1}{d(\alpha)}\log M(\alpha)$ sa hauteur logarithmique absolue. 

Pour des nombres $R \in \reels \cup \{\infty\}$ et $r \in \reels$ tels que $R>r>0,$ et une fonction $f$ analytique sur $\DRouvert,$ on note, pour tout entier $D \geq 1$ et tout r\'eel $N \geq 0, \esp \Sigma_{D,N}=\Sigma _{D,N}(f,r)$ l'ensemble des nombres $\alpha \in \Qbarre \cap  \Drferme$ tels que 
$$f(\alpha) \in
\Qbarre ,\esp [\rationnels(\alpha , f(\alpha )):\rationnels] \leq D ,\hspace{0,2cm} h(\alpha) \leq N \esp \mathrm{et} \esp h(f(\alpha)) \leq N,$$
et $\sigma_{D,N}$ son cardinal. Ainsi, $S_f$ est la r\'eunion des $\Sigma_{D,N}$ pour $D\geq 1$ et $N\geq 0.$

\begin{theoreme}\label{contrex}
Soit $\phi$ une fonction positive telle que $\phi(x)/x$ tende vers $0$ quand $x \to \infty.$ Il existe une fonction $f$ enti\`ere et transcendante sur $\complexes [z ],$ v\'erifiant
$$f(\alpha) \in \rationnels[\alpha], \esp \forall \alpha \in \Qbarre, $$
et telle que, pour tout entier $D \geq 1,$ il existe une infinit\'e de r\'eels $N \geq 0,$ v\'erifiant
$$\card(\Sigma _{D,N}(f,1)) \geq \frac{1}{2}e^{D(D+1)\phi(N)}.$$

\end{theoreme}

Pour $D$ et $N$ fix\'es, cet ensemble est fini. En effet, il est contenu dans 
%
l'ensemble
$$E_{D,N}= \{ \alpha \in \Qbarre / \hspace{0,2cm}[\rationnels (\alpha
):\rationnels ] \leq D , \esp h(\alpha )\leq N \}$$
qui est fini. Notons $\epsilon_{D,N}$ son cardinal. 

\medskip 

\begin{lemme}\label{e_{D,N}}
Pour tout entier $D \geq 1$ et tout r\'eel $N \geq 0,$ le cardinal $\epsilon _{D,N}$ de $E_{D,N}$ v\'erifie
$$e^{D(D+1)(N-1)} < \epsilon _{D,N} \leq e^{D(D+1)(N+1)}.$$

\end{lemme} 

La d\'emonstration est \'el\'ementaire. La majoration est facile ; une estimation plus pr\'ecise quand $N$ tend vers l'infini est connue \cite{chern-vaaler}, mais elle n'est pas explicite. Pour la minoration, on compte les polyn\^omes 
dans $\enteros[X]$ de hauteur born\'ee qui sont 2-Eisenstein sur l'anneau $\enteros_2$ des entiers 2-adiques ; une minoration l\'eg\'erement plus pr\'ecise se trouve aussi dans \cite{loher}. Aucune estimation asymptotique ne semble actuellement connue (\cite{schmidt} p. 27).

Le r\'esultat suivant constitue le th\'eor\`eme principal de ce travail.

\begin{theoreme}{Majoration du cardinal de $\Sigma_{D,N}(f,r).$}\label{propprincipale}

Soient $R \in \reels \cup \{\infty\}$ et $ r \in \reels$ tels que $R>r>0,$ et $f$ une fonction analytique sur $\DRouvert.$ Il existe une constante $\gamma >0,$ d\'ependant uniquement de $R, r$ et $f,$ telle que, pour tout entier $D \geq1,$ il existe une infinit\'e de r\'eels $N \geq 0$ pour lesquels
$$\card(\Sigma _{D,N}(f,r)) < \gamma D^3N^2.$$

\end{theoreme}

Le Th\'eor\`eme \ref{contrex} (resp. \ref{propprincipale})  montre qu'on ne peut pas remplacer, dans le Th\'eor\`eme \ref{propprincipale} (resp. \ref{contrex}), ``il existe une infinit\'e de r\'eels $N$'' par ``pour tout $N$ assez grand'', et ceci m\^eme pour une borne sup\'erieure allant jusqu'\`a $\frac{1}{2}e^{D(D+1)\phi(N)}.$ Pour le cas particulier o\`u $D=1,$ Elkies (\cite{elkies} th.4.) a d\'emontr\'e que, pour tout $\epsilon > 0,$ il existe une constante $A_{\epsilon} > 0,$ telle que, pour tout r\'eel $N \geq 0, \esp \card(\Sigma _{1,N}) \leq A_{\epsilon} e^{\epsilon N}.$

\section{D\'emonstration du Th\'eor\`eme \ref{contrex}}


Soient $\phi$ une fonction v\'erifiant les hypoth\`eses du Th\'eor\`eme \ref{contrex}, $(b_k)_{k\geq 1}$ une suite de nombres r\'eels $>0$ et telle que la s\'erie $\sum_{k\geq1}b_k$ soit convergente et $x_0 \geq 1$ un nombre r\'eel tel que pour tout r\'eel $x \geq x_0,$
\begin{equation}
\phi(x) \leq x-1. \label{x_0}
\end{equation}

Soient $(N_{\delta})_{\delta\geq 1}$ une suite strictement croissante de nombres r\'eels $\geq x_0$ tendant vers l'infini et $(a_{\delta})_{\delta \geq 1}$ une suite de nombres rationnels v\'erifiant les conditions suivantes.

\bigskip

\noindent i)  Pour une infinit\'e de $\delta, \esp a_{\delta} \ne 0$ et pour tout $\delta \geq 1, \esp |a_{\delta}| \leq b_{\delta}(\delta + e^{\delta N_{\delta}})^{-\delta \epsilon_{\delta,N_{\delta}}}.$

\noindent ii) Pour tout $\delta \geq 2,$ 
\begin{equation}
N_{\delta} \geq 2 \left( \log(\delta-1) + \sum_{k=1}^{\delta-1} h(a_k) + (\delta-1)^2 \epsilon_{\delta-1,N_{\delta-1}}(\log2 + 1 + N_{\delta-1}) \right).  \label{Ndelta>.} 
\end{equation}

\noindent iii) Pour tout $\delta \geq 2,$
\begin{equation}
\frac{N_{\delta}}{\phi(N_{\delta})} \geq 2 (\delta-1)^2 \epsilon_{\delta-1,N_{\delta-1}}.                                                   \label{Ndelta>2 Phi}
\end{equation}

On pose, pour tout $k \geq 1,$ et tout $z\in \complexes,$ 
$$P_k(X) = \prod_{\beta \in E_{k,N_k}}(X - \beta ) \esp \mathrm{et} \esp f(z) = \sum_{k\geq 1}a_k P_k(z).$$

\medskip

De la condition i) on d\'eduit que la fonction $f$ est une fonction enti\`ere et non polyn\^omiale; elle est donc  transcendante.

\medskip


Soit $\alpha \in \Qbarre.$ Montrons que $f(\alpha) \in \rationnels (\alpha).$ Par hypoth\`ese, la suite $(N_{\delta})_{\delta \geq 1}$ tend vers l'infini, donc les ensembles $(E_{k,N_k})_{k\geq 1}$ forment une suite croissante pour l'inclusion, de r\'eunion $\Qbarre,$ et il existe $k_0 \geq 1$ tel que, pour tout $k \geq k_0, P_k(\alpha) = 0$ et 
$$f(\alpha) = \sum_{k=1}^{k_0 - 1}a_k P_k(\alpha).$$

De plus, les polyn\^omes $P_k$ sont fix\'es par tout \'el\'ement du groupe de Galois de $\rationnels(\alpha)$ sur $\rationnels;$ ils sont donc \`a coefficients rationnels, et comme, par d\'efinition, les nombres $a_k$ sont aussi rationnels, $f(\alpha) \in \rationnels[\alpha].$


Soient $D$ et $d$ deux entiers tels que $d \geq D \geq 1.$ Montrons que l'ensemble
$$\Sigma_{D,N_d}(f,1)= \{\alpha \in \Qbarre \cap \Duferme / \esp \deg(\alpha) \leq D, \esp h(\alpha) \leq N_d, \esp h(f(\alpha)) \leq N_d \}$$
contient $E_{D,\phi(N_d)+1} \cap \Duferme.$ Soit $\alpha \in \Qbarre$ tel que 
$$\deg(\alpha) \leq d \esp \textrm{et} \esp  h(\alpha) \leq \phi(N_d) + 1.$$

Comme $N_d \geq x_0,$ alors, d'apr\`es (\ref{x_0}), $\phi(N_d) \leq N_d - 1,$ et $h(\alpha) \leq N_d.$
Comme, en plus, $\deg (\alpha) \leq d,$ alors $\alpha \in E_{d, N_d}$ et $\forall k \geq d, \hspace{0,2cm} P_k(\alpha) = 0.$

Si $d=D=1,$ alors $f(\alpha)=0$ et on a l'inclusion $E_{1,\phi(N_1)+1}\cap \Duferme \subset \Sigma _{1,N_1}(f,1).$

Supposons maintenant que $d>D \geq 1.$ Alors 
$$f(\alpha) = \sum_{k=1}^{d-1} a_k P_k(\alpha).$$

D'apr\`es les conditions (\ref{Ndelta>.}) et (\ref{Ndelta>2 Phi}), on montre que la hauteur de $f(\alpha)$ est inf\'erieure \`a 
$$\phi(N_d)\left( (d-1)^2 \epsilon_{d-1, N_{d-1}} \right) + \log(d-1) + A_{d-1} + (d-1)^2 \epsilon_{d-1, N_{d-1}} (\log 2 + 1 + N_{d-1}) \leq N_d.$$ 

En particulier, si $\alpha \in E_{D,\phi(N_d)+1} \cap \Duferme ,$ alors $\deg(\alpha) \leq D < d$ et $h(f(\alpha)) \leq N_d,$ ce qui montre que $\alpha \in \Sigma _{D, N_d}(f,1).$ 

Ainsi, pour $d \geq D \geq 1,$ nous avons
$$\sigma_{D,N_d} = \card(\Sigma_{D,N_d} (f,1)) \geq \card(E_{D,\phi(N_d)+1} \cap \Duferme).$$

En remarquant que si $\alpha \in \Qbarre-\{0\}$ on a $d(\alpha) = d(1/\alpha)$ et $h(\alpha) = h(1/\alpha),$ et en appliquant le lemme \ref{e_{D,N}}, nous obtenons
$$\sigma_{D,N_d} \geq \frac{1}{2} \card(E_{D,\phi(N_d)+1}) > \frac{1}{2} e^{D(D+1)\phi(N_d)}.$$
%


\section{D\'emonstration du Th\'eor\`eme \ref{propprincipale}}


\subsection{\textbf{Choix des param\`etres et construction de
  $S_{D,N}$}}\label{choix des param}


On se donne $R, r$ et $f$ v\'erifiant les hypoth\`eses du th\'eor\`eme. On pose $c_0=\log \left( \frac{R^2 + r^2}{2rR} \right).$ C'est un nombre r\'eel strictement positif. On se donne aussi un entier $\gamma$ strictement positif et v\'erifiant $\gamma > 2\left(\frac{6}{c_0} \right)^{2}.$

Par l'absurde, on suppose qu'il existe un entier $D \geq 1$ et un entier $N_0 \geq 1$ suffisamment grand tels que pour tout $N \geq N_0$ on ait $\card(\Sigma _{D,N}(f,r)) \geq \gamma D^3N^2.$

Pour tout $N\geq N_0-1,$ on extrait  de $\Sigma_{D,N}$ un sous-ensemble $S_{D,N}$ dont le nombre d'\'el\'ements, que l'on note $s_N,$ est exactement $\gamma D^3N^2.$ On pose 
$$T = \left[ \frac{c_0 \gamma}{6} D^2 N_0 \right] \esp \mathrm{et} \esp u_1 = \log\left( 2T^4 e^{2N_0T} |f|_R^T\right).$$

Comme $N_0$ est suffisamment grand on a 
\begin{equation}
c_0 s_{N_0-1} > u_1 + (D-1)\log(2T^4) + 2TN_0(D-1) + 2TN_0D.  \label{N_0}
\end{equation}


\subsection{\textbf{La fonction auxiliaire }}

Un lemme de Siegel \cite{gmw} nous donne l'existence d'un polyn\^ome $P \in \enteros[X,Y],$ non nul (et d\'ependant de $N_0$), dont le degr\'e par rapport \`a chacune des variables est $\leq T,$ tel que 
$$P(\alpha, f(\alpha))=0,\hspace{0,5cm}  \forall \alpha \in S_{D,N_0},$$
et dont les coefficients sont major\'es en valeur absolue par $2T^2e^{2TN_0}.$

La fonction $F(z)=P(z,f(z))$ s'annule sur tout l'ensemble $S_{D,N_0}.$


\subsection{\textbf{R\'ecurrence et conclusion}}

On montre par r\'ecurrence que pour tout $N \geq N_0, \esp F$ s'annule sur tout l'ensemble $S_{D,N}.$ 

On se donne $N \geq N_0$ et on suppose que $F(\alpha) = 0,$ pour tout $\alpha \in S_{D,N}.$ En appliquant un lemme de Schwarz \cite{w}, et gr\^ace \`a la borne obtenue pour les coefficients de $P,$ on montre que, pour tout $\alpha \in \Drferme,$ et donc, en particulier, pour tout $\alpha \in S_{D,N+1},$ 
$$|F(\alpha)| \leq e^{u_1 - c_0s_N}.$$

Gr\^ace \`a la majoration de la hauteur des coefficients de $P,$ et de l'in\'egalit\'e (\ref{N_0}), et, en notant $L(P)$ la longueur de $P,$ on obtient
$$|F(\alpha)| < L(P)^{-(D-1)} e^{-DT(h(\alpha)+h(f(\alpha))}.$$
L'in\'egalit\'e de Liouville \cite{w} implique $F(\alpha)=0.$

On a ainsi montr\'e que  pour tout $N \geq N_0,$ la fonction $F$ s'annule sur $S_{D,N},$ ce qui contredit le th\'eor\`eme sur les z\'eros isol\'es d'une fonction holomorphe et nous donne le r\'esultat cherch\'e.


%

\medskip

\noindent Andrea Surroca\\
Institut de Math\'ematiques de Jussieu \\ 
\'Equipe de Th\'eorie des nombres, Case 247 \\
175, rue du Chevaleret 75013 Paris.\\
surroca@math.jussieu.fr

\end{document}